\documentclass[a4paper,reqno,11pt]{amsart}

\usepackage[ansinew]{inputenc}
\usepackage{amsfonts}
\usepackage{amsmath}
\usepackage{amssymb}
\usepackage{latexsym}

\newtheorem{theorem}{Theorem}

\newtheorem{corollary}[theorem]{Corollary}

\theoremstyle{definition}

\numberwithin{equation}{section}

\newcommand{\D}{\mathbb{D}}
\newcommand{\C}{\mathbb{C}}
\newcommand{\NN}{\mathbb{N}}
\newcommand{\N}{\mathcal{N}}
\newcommand{\M}{\mathcal{M}}

\newcommand{\DD}{\mathcal{D}}
\newcommand{\vp}{\varphi}
\newcommand{\om}{\omega}

\def\UBC{\mathord{\rm UBC}}

\renewcommand{\H}{\mathcal{H}}

\theoremstyle{remark}

\begin{document}
\title[Growth estimates for solutions of algebraic differential equations]{Growth estimates for meromorphic solutions of higher order algebraic differential equations}

\author{Shamil Makhmutov}
\address{Department of Mathematics, College of Science, Sultan Qaboos University, P.O. Box 36, PC 123 Al Khodh, Muscat, Sultanate of Oman}
\email{makhm@squ.edu.om}

\author{Jouni R\"atty\"a}
\address{University of Eastern Finland, P.O.Box 111, 80101 Joensuu, Finland}
\email{jouni.rattya@uef.fi}

\author{Toni Vesikko}
\address{University of Eastern Finland, P.O.Box 111, 80101 Joensuu, Finland}
\email{tonive@uef.fi}

\maketitle

\begin{abstract}
We establish pointwise growth estimates for the spherical derivative of solutions of the first order algebraic differential equations. A generalization of this result to higher order equations is also given. We discuss the related question of when for a given class $X$ of meromorphic functions in the unit disc, defined by means of the spherical derivative,  and $m\in\NN\setminus\{1\}$, $f^m\in X$ implies $f\in X$. An affirmative answer to this is given for example in the case of $\UBC$, the $\alpha$-normal functions with $\alpha\ge1$ and certain (sufficiently large) Dirichlet type classes. 
\end{abstract}

\section{Introduction and main results}

Let $\H(\D)$ and $\M(\D)$ denote the sets of analytic and meromorphic functions in the unit disc $\D=\{z\in\C:|z|<1\}$, respectively. For $n,N\in\NN$, consider the $N$-th order algebraic differential equation
	\begin{equation}\label{orderNref}
	\left(f^{(N)}\right)^n + \sum_{k=1}^{n} P_{k,N}\left(f\right)\left(f^{(N)}\right)^{n-k} = 0,
	\end{equation}
where 
	\begin{equation}
	P_{k,N}(f) 
	=\sum_{j_0=0}^{m_{k,0}}\sum_{j_1=0}^{m_{k,1}}\cdots\sum_{j_{N-1}=0}^{m_{k,N-1}} 
	a_{k,j_0,\ldots,j_{N-1}}\prod_{\ell=0}^{N-1} \left(f^{(\ell)}\right)^{j_\ell},\quad k=1,\ldots,n,
	\end{equation}
with $a_{k,j_0,\ldots,j_{N-1}}\in\H(\D)$ and $m_{k,j}\in\NN\cup\{0\}$ for all $j=0,\ldots,N-1$ and $k=1,\ldots,n$. The case $N=1$ reduces to the first order equation 
    \begin{equation}\label{eq-1}
     (f')^n + \sum_{k=1}^{n} P_k(f)\,(f')^{n-k} = 0,
    \end{equation}
where
    \begin{equation}\label{Pk-1}
    P_k(f)=\sum_{j=0}^{m_k} a_{k,j} f^j , \quad 1\le k \le n\in\NN,\quad m_k\in\NN\cup\{0\},
    \end{equation}
and $a_{k,j}\in\H(\D)$ for all $k=1,\ldots,n$ and $j=0,\ldots,m_k$. 

The main result of this study is a pointwise growth estimate for the spherical derivative of meromorphic solutions of \eqref{orderNref}. The method of proof does not depend on the underlying domain and can be performed on any set. The normal and Yosida solutions of algebraic differential equations similar to \eqref{orderNref} and \eqref{eq-1} have been studied extensively via techniques like, for example, the Lohwater-Pommerenke method \cite{A-M-R2010,AW2011}. There are multiple existing results concerning normality conditions and the behaviour of the spherical derivatives of the solutions \cite{A-M-R2010,L-L-Y2003,Makhmutov2011}. The method we employ allows us to consider solutions in classes which are strictly smaller than the class of normal functions.

Before stating the results, a word about the notation used. The letter $C=C(\cdot)$ will denote an absolute constant whose value depends on the parameters indicated
in the parenthesis, and may change from one occurrence to another.
We will use the notation $a\lesssim b$ if there exists a constant
$C=C(\cdot)>0$ such that $a\le Cb$, and $a\gtrsim b$ is understood
in an analogous manner. In particular, if $a\lesssim b$ and
$a\gtrsim b$, then we write $a\asymp b$ and say that $a$ and $b$ are comparable.

\begin{theorem}\label{thm:main}
Let $n,N\in\NN$ and $M_\ell\in\NN\cup\{0\}$ such that 
	\begin{equation}\label{hypotheses}
	M_0\ge\max_{k=1,\ldots,n}\frac{m_{k,0}}{k}-2\quad\textrm{and}\quad
	M_\ell\ge\max_{k=1,\ldots,n}\frac{m_{k,\ell}}{k}-1,\quad \ell=1,\ldots,N-1.
	\end{equation}
Then each meromorphic solution $f$ of \eqref{orderNref} satisfies
	\begin{equation*}
	\prod_{\ell=0}^{N-1}\left(\left(f^{(\ell)}\right)^{M_{\ell}+1}\right)^{\#} 
	\lesssim\sum_{k=1}^n\left(\sum_{j_0=0}^{m_{k,0}}\sum_{j_1=0}^{m_{k,1}}\cdots\sum_{j_{N-1}=0}^{m_{k,N-1}}
	|a_{k,j_0,\cdots,j_{N-1}}|^{\frac{1}{k}}\right).
	\end{equation*}
In particular, if $\max_{k=1,\ldots,n}\frac{m_k}{k}\le M_0+2$, then each meromorphic solution $f$ of \eqref{eq-1} satisfies
	\begin{equation*}
	\left(f^{M_0+1}\right)^{\#} 
	\lesssim\sum_{k=1}^n\sum_{j=0}^{m_{k}}|a_{k,j}|^{\frac{1}{k}}.
	\end{equation*}
\end{theorem}

Theorem~\ref{thm:main} has the following immediate consequence which deserves to be stated separately.

\begin{corollary}\label{corollary2}
Let $n,N\in\NN$, $\max_{k=1,\ldots,n}\frac{m_{k,0}}{k}\le2$ and $\max_{k=1,\ldots,n}\frac{m_{k,\ell}}{k}\le1$ for all $\ell=1,\ldots,N-1$. Then each meromorphic solution $f$ of \eqref{orderNref} satisfies
	\begin{equation}\label{Eq:corollary}
	\prod_{\ell=0}^{N-1}\left(f^{(\ell)}\right)^{\#} 
	\lesssim\sum_{k=1}^n\left(\sum_{j_0=0}^{m_{k,0}}\sum_{j_1=0}^{m_{k,1}}\cdots\sum_{j_{N-1}=0}^{m_{k,N-1}}
	|a_{k,j_0,\cdots,j_{N-1}}|^{\frac{1}{k}}\right).
	\end{equation}
\end{corollary}

The product of the spherical derivatives on the left hand side of \eqref{Eq:corollary} appears in a natural way in the study of (weighted) normal functions \cite{CL1996,CL1998,Grohn2017,Lappan-77,Makhmutov1986,Xu2000}. 

The second part of Theorem~\ref{thm:main} gives arise to the question of when for a given class $X\subset\M(\D)$ and $m\in\NN\setminus\{1\}$, $f^m\in X$ implies $f\in X$. An immediate observation is that, roughly speaking, this implication cannot be true if $X$ is sufficiently small and defined in terms of the spherical derivative. More precisely, for $f_p(z)=(1-z)^{-p}$ with $0<p<1/m<1$ we have $f_p^\#(z)\asymp|1-z|^{p-1}$ and $(f_p^m)^\#\asymp|1-z|^{mp-1}$ as $\D\ni z\to1$. This shows that $(f_p^m)^\#$ is essentially smaller than $(f_p)^\#$ when $z\to1$, yet of course $f_p\in\N^{1-p}\subset\N$ for all $0<p<1$. Recall that for $0<\alpha<\infty$, the class $\N^\alpha$ of $\alpha$-normal functions is the set of functions $f\in\M(\D)$ such that
	$$
	\|f\|_{\N^\alpha}=\sup_{z\in\D}f^\#(z)(1-|z|^2)^\alpha<\infty,
	$$
and its subset $\N^\alpha_0$ of strongly $\alpha$-normal functions consists of functions $f\in\M(\D)$ such that $f^\#(z)(1-|z|^2)^\alpha\to0$
 as $|z|\to1^-$.

We next offer an affirmative answer to the question of when $f^m\in X$ implies $f\in X$ in the case of certain function classes. To do this, definitions are needed. An increasing function $\vp:[0,1)\to(0,\infty)$ is smoothly increasing if
$\vp(r)(1-r)\to\infty$, as $r\to1^-$, and
		\begin{equation*}
		\frac{\vp(|a+z/\vp(|a|)|)}{\vp(|a|)}\to1,\quad |a|\to1^-,
    \end{equation*}
uniformly on compact subsets of $\C$. For such a $\vp$, a
function $f\in\M(\D)$ is $\vp$-normal if
    \begin{equation}\label{definition}
    \|f\|_{\N^\vp}=\sup_{z\in\D}\frac{f^{\#}(z)}{\vp(|z|)}<\infty,
    \end{equation}
and strongly $\vp$-normal if $\frac{f^{\#}(z)}{\vp(|z|)}\to0$, as $|z|\to1^-$. The classes of $\vp$-normal and strongly $\vp$-normal functions are denoted by $\N^\vp$ and $\N^\vp_0$, respectively.

Let $\om\in L^1(0,1)$. The extension defined by $\om(z)=\om(|z|)$ for all $z\in\D$ is called a radial weight on $\D$. For such an $\om$, denote 
	$$
	\om^\star(z)=\int_{|z|}^1\om(s)\log\frac{s}{|z|}s\,ds,\quad z\in\D\setminus\{0\}.
	$$	
The Dirichlet class $\DD^\#_{\om^\star}$ consists of $f\in\M(\D)$ such that 
	$$
	\int_\D f^\#(z)^2\om^\star(z)\,dA(z)<\infty,
	$$
where $dA(z)=rdrd\theta/\pi$ for $z=re^{i\theta}$. Moreover, the Dirichlet class $\DD^\#_{\alpha}$ consists of $f\in\M(\D)$ such that 
	$$
	\int_\D f^\#(z)^2(1-|z|^2)^\alpha\,dA(z)<\infty.
	$$
A function $f\in\M(\D)$ belongs to $\UBC$ if 
	$$
	\sup_{a\in\D}\int_\D f^\#(z)^2\log\frac1{|\varphi_a(z)|}\,dA(z)<\infty,
	$$
where $\vp_a(z)=(a-z)/(1-\overline{a}z)$.

With these preparations we can state our next result.

\begin{theorem}\label{theorem:f-m-implications}
Let $m\in\NN\setminus\{1\}$, $1<p<\infty$, $\vp$ a smoothly increasing and $\om$ a radial weight. Then the following statements are valid:
	\begin{itemize}
	\item[\rm(i)] $f^m\in\N\Rightarrow f\in\N$;
	\item[\rm(ii)] $f^m\in\N_0\Rightarrow f\in\N_0$;
	\item[\rm(iii)] $f^m\in\N^\vp\Rightarrow f\in\N^\vp$;
	\item[\rm(iv)] $f^m\in\N^\vp_0\Rightarrow f\in\N^\vp_0$;
	\item[\rm(v)] $f^m\in\DD^\#_1\Rightarrow f\in\DD^\#_1$;
	\item[\rm(vi)] $f^m\in\DD^\#_{\om^\star}\Rightarrow f\in\DD^\#_{\om^\star}$;
	\item[\rm(vii)] $f^m\in\UBC\Rightarrow f\in\UBC$.
	\end{itemize}
\end{theorem}

The proofs of the cases (i)--(iv) are based on the so-called five-point theorems, named by the celebrated result for normal functions due to Lappan \cite{Lappan1974}. Such results exist also in the setting of meromorphic functions on the whole plane and are usually given in terms of Yosida and $\vp$-Yosida functions, see \cite{A-M-R2010,Makhmutov2001}. Therefore we may obtain analogues of the statements (i)-(iv) for those classes; details are left for interested reader. 

The argument of proof used in (v)--(vii) uses \cite[Theorem~2]{Yamashita1981} due to Yamashita on meromorphic Hardy classes. This is not exclusive for the disc either and can be performed also on the plane. The argument yields the inequality
	$$
	\int_\C f^\#(z)^2\left(\int_{|z|}^\infty\log\frac{r}{|z|}\om(r)\,dr\right)dA(z)
	\lesssim\int_\C\left(f^m\right)^\#(z)^2\left(\int_{|z|}^\infty\log\frac{r}{|z|}\om(r)\,dr\right)dA(z)+1,
	$$
valid for all meromorphic functions $f$ and radial weights $\om$ on $\C$. This is an analogue of (vi) for $\C$. One natural choice for $\om$ in this case is $\om(z)=(1+|z|)^{-\alpha}$ for $1<\alpha<\infty$.

Let $0<\alpha<1$ be fixed. By Theorem~\ref{theorem:f-m-implications}(i) we know that if $f^m\in\N^\alpha\subset\N$, then $f\in\N$. For $0<\alpha<1$ and $m\in\NN\setminus\{1\}$, write 
	$$
	\beta_{\alpha,m}=\inf_{0<\gamma\le1}\left\{f\in\N^\gamma:f^m\in\N^\alpha\right\}.
	$$
The function $f_{p}$ with $p=\frac{1-\alpha}{m}$ considered after Corollary~\ref{corollary2} and Theorem~\ref{theorem:f-m-implications}(i) show that $1-\frac{1-\alpha}{m}\le\beta_{\alpha,m}\le1$. The exact value of $\beta_{\alpha,m}$ is unknown. 

Theorem~\ref{theorem:f-m-implications} shows that $f^m\in\DD^\#_\alpha$ with $\alpha\ge1$ implies $f\in\DD^\#_\alpha$. Further, the function $f_p$ with $-\frac{\alpha}{2m}<p\le-\frac{\alpha}{2}$ shows that $f^m\in\DD^\#_\alpha$ with $\alpha<0$ does not imply $f\in\DD^\#_\alpha$. It is natural to ask what happens with the range $0\le\alpha<1$? We do not know an answer to this question. 

We next aim for combining Theorems~\ref{thm:main} and \ref{theorem:f-m-implications} in order to find a set of sufficient conditions for the coefficients of \eqref{orderNref} that force meromorphic solution $f$ to belong to certain function classes. To do this, some more notation is needed. For $0<p<\infty$, the weighted growth space $H^\infty_p$ consists of $f\in\H(\D)$ such that
		$$
		\|f\|_{H^\infty_p}=\sup_{z\in\D}|f(z)|(1-|z|)^p<\infty.
		$$
Similarly, $H^p_\vp$ consists of $f\in\H(\D)$ such that
		$$
		\|f\|_{H^p_\vp}=\sup_{z\in\D}\frac{|f(z)|}{\vp(|z|)}<\infty.
		$$

For $0<p<\infty$ and a radial weight $\om$, the Bergman space $A^p_\om$ consists of $f\in\H(\D)$ such that
		$$
		\|f\|_{A^p_\om}^p=\int_\D|f(z)|^p\om(z)\,dA(z)<\infty.
		$$
Note that the Hardy-Spencer-Stein formula yields
	$$
	\|f\|_{A^p_\omega}^p=p^2\int_{\D}|f(z)|^{p-2}|f'(z)|^2\omega^\star(z)\,dA(z)+\omega(\D)|f(0)|^p,\quad f\in\H(\D),
	$$
by \cite[Theorem~4.2]{PelaezRattya2014}. This explains how the associated weight $\om^\star$ raises in a natural manner.

The following result is an immediate consequence of Theorems~\ref{thm:main} and \ref{theorem:f-m-implications}. 

\begin{corollary}
Let $n\in\NN$ and $M_0\in\NN\cup\{0\}$ such that $\max_{k=1,\ldots,n}\frac{m_k}{k}\le M_0+2$, and let $\vp$ be smoothly increasing. 
	\begin{itemize}
	\item[\rm(i)] If $a_{k,j}\in H^\infty_k$ (resp. $a_{k,j}\in H^\infty_{k,0}$) for all $j=0,\ldots,m_k$ and $k=1,\ldots,n$, then each meromorphic solution $f$ of \eqref{eq-1} belongs to $\N$ (resp. $\N_0)$.
	\item[\rm(ii)] If $a_{k,j}\in H^\infty_{\vp^k}$ (resp. $a_{k,j}\in H^\infty_{\vp^k,0}$) for all $j=0,\ldots,m_k$ and $k=1,\ldots,n$, then each meromorphic solution $f$ of \eqref{eq-1} belongs to $\N^\vp$ (resp. $\N^\vp_0$).
	\item[\rm(iii)] If $a_{k,j}\in A^{\frac2k}_{\om^*}$ for all $j=0,\ldots,m_k$ and $k=1,\ldots,n$, then each meromorphic solution $f$ of \eqref{eq-1} belongs to $\DD^\#_{\om^*}$.
		\item[\rm(iv)] If 
			$$
			\sup_{a\in\D}\int_\D|a_{k,j}(z)|^\frac{2}{k}\log\frac{1}{|\vp_a(z)|}\,dA(z)<\infty
			$$
	for all $j=0,\ldots,m_k$ and $k=1,\ldots,n$, then each meromorphic solution $f$ of \eqref{eq-1} belongs to $\UBC$.
	\end{itemize}
\end{corollary}

It is well known that in (iv) one may replace the Green's function $\log\frac{1}{|\vp_a(z)|}$ by the term $1-|\vp_a(z)|^2$ in the statement. This is due to the analyticity of the coefficients.

\section{Proof of Theorem~\ref{thm:main}}

We first multiply \eqref{orderNref} by 
	$$
	\left(\left(M_{N-1}+1\right)\left(f^{N-1}\right)^{M_{N-1}}\prod_{\ell=0}^{N-2}\left(\left(f^{(\ell)}\right)^{M_\ell + 1}\right)'\right)^n
	$$ 
to obtain
	\begin{equation*}
	\begin{split} 
	\prod_{\ell=0}^{N-1}\left(\left(\left( f^{(\ell)}\right)^{M_{\ell}+1}\right)'\right)^n
	+&\sum_{k=1}^nP_{k,N}(f)\left(f^{(N)}\right)^{n-k}\left(M_{N-1}+1\right)^n\left(f^{(N-1)}\right)^{M_{N-1}n}\\
	&\cdot\left(\prod_{\ell=0}^{N-2}\left(\left(f^{(\ell)}\right)^{M_\ell+1}\right)'\right)^n=0,
	\end{split}
	\end{equation*}
and then divide it by $\prod_{\ell=0}^{N-1}\left(1+|f^{(\ell)}|^{2\left(M_\ell+1\right)}\right)^n$ to get
	\begin{equation*}
	\begin{split} 
	\prod_{\ell=0}^{N-1}\left(\frac{\left(\left( f^{(\ell)}\right)^{M_{\ell}+1}\right)'}{1+|f^{(\ell)}|^{2\left(M_\ell+1\right)}}\right)^n
	+&\sum_{k=1}^nP_{k,N}(f)\left(M_{N-1}+1\right)^k\left(f^{(N-1)}\right)^{M_{N-1}k}\\
	&\cdot\frac{\left(\prod_{\ell=0}^{N-2}\left(\left(f^{(\ell)}\right)^{M_\ell+1}\right)'\right)^n
	\left(\left(\left(f^{(N-1)}\right)^{M_{N-1}+1}\right)'\right)^{n-k}}
	{\prod_{\ell=0}^{N-1}\left(1+|f^{(\ell)}|^{2\left(M_\ell+1\right)}\right)^n}=0.
	\end{split}
	\end{equation*}
By reorganizing terms and taking moduli, we deduce  
	\begin{equation*}
	\begin{split}
	\prod_{\ell=0}^{N-1}\left(\left(f^{(\ell)}\right)^{M_{\ell}+1}\right)^\#
	&\le\Bigg(\sum_{k=1}^n\left|P_{k,N}(f)\right|\left(M_{N-1}+1\right)^k\left|f^{(N-1)}\right|^{M_{N-1}k}\\
	&\quad\cdot\frac{\prod_{\ell=0}^{N-2}\left(M_\ell+1\right)^n|f^{(\ell)}|^{M_\ell n}|f^{(\ell+1)}|^n}{\prod_{\ell=0}^{N-1}\left(1+|f^{(\ell)}|^{2\left(M_\ell+1\right)}\right)^n}\left|\left(\left(f^{(N-1)}\right)^{M_{N-1}+1}\right)'\right|^{n-k}\Bigg)^\frac1n\\
	&\le\sum_{k=1}^n\left|P_{k,N}(f)\right|^\frac1n\left(M_{N-1}+1\right)^\frac{k}n\left|f^{(N-1)}\right|^{\frac{M_{N-1}k}{n}}\\
	&\quad\cdot\frac{\prod_{\ell=0}^{N-2}\left(M_\ell+1\right)|f^{(\ell)}|^{M_\ell}|f^{(\ell+1)}|}
	{\prod_{\ell=0}^{N-1}\left(1+|f^{(\ell)}|^{2\left(M_\ell+1\right)}\right)}
	\left|\left(\left(f^{(N-1)}\right)^{M_{N-1}+1}\right)'\right|^{\frac{n-k}{n}}\\
	&=\sum_{k=1}^n\left|P_{k,N}(f)\right|^\frac1n\left(M_{N-1}+1\right)^\frac{k}n
	\prod_{\ell=0}^{N-2}\left(\left(\left(f^{(\ell)}\right)^{M_\ell+1}\right)^\#\right)^\frac{k}{n}\\
	&\quad\cdot\left(\frac{|f^{(N-1)}|^{M_{N-1}}}{1+|f^{(N-1)}|^{2(M_{N-1}+1)}}\right)^{\frac{k}{n}}
	\left(\prod_{\ell=0}^{N-1}\left(\left(f^{(\ell)}\right)^{M_{\ell}+1}\right)^\#\right)^{\frac{n-k}{n}}.
	\end{split}
	\end{equation*}
It follows that
	\begin{equation*}
	\begin{split}
	\prod_{\ell=0}^{N-1}\left(\left(f^{(\ell)}\right)^{M_{\ell}+1}\right)^\#
	&\le\Bigg(\sum_{k=1}^n\left|P_{k,N}(f)\right|^\frac1n\left(M_{N-1}+1\right)^\frac{k}n
	\prod_{\ell=0}^{N-2}\left(\left(\left(f^{(\ell)}\right)^{M_\ell+1}\right)^\#\right)^\frac{k}{n}\\
	&\quad\cdot\left(\frac{|f^{(N-1)}|^{M_{N-1}}}{1+|f^{(N-1)}|^{2(M_{N-1}+1)}}\right)^{\frac{k}{n}}
	\Bigg)^{\frac{n}{k}}\\
	&\lesssim\sum_{k=1}^n\left|P_{k,N}(f)\right|^\frac1k\left(M_{N-1}+1\right)
	\prod_{\ell=0}^{N-2}\left(\left(f^{(\ell)}\right)^{M_\ell+1}\right)^\#\\
	&\quad\cdot\frac{|f^{(N-1)}|^{M_{N-1}}}{1+|f^{(N-1)}|^{2(M_{N-1}+1)}},
	\end{split}
	\end{equation*}
where
	\begin{equation*}
	\begin{split}
	\left|P_{k,N}(f)\right|^\frac1k
	&\le\left(\sum_{j_0=0}^{m_{k,0}}\sum_{j_1=0}^{m_{k,1}}\cdots\sum_{j_{N-1}=0}^{m_{k,N-1}} 
	\left|a_{k,j_0,\ldots,j_{N-1}}\right|\prod_{\ell=0}^{N-1}\left|f^{(\ell)}\right|^{j_\ell}\right)^\frac1k\\
	&\le\sum_{j_0=0}^{m_{k,0}}\sum_{j_1=0}^{m_{k,1}}\cdots\sum_{j_{N-1}=0}^{m_{k,N-1}} 
	\left|a_{k,j_0,\ldots,j_{N-1}}\right|^\frac1k\prod_{\ell=0}^{N-1}\left|f^{(\ell)}\right|^{\frac{j_\ell}{k}}.
	\end{split}
	\end{equation*}
Hence
	\begin{equation*}
	\begin{split}
	\prod_{\ell=0}^{N-1}\left(\left(f^{(\ell)}\right)^{M_{\ell}+1}\right)^\#
	&\lesssim\left(M_{N-1}+1\right)\sum_{k=1}^n\left(\sum_{j_0=0}^{m_{k,0}}\sum_{j_1=0}^{m_{k,1}}\cdots\sum_{j_{N-1}=0}^{m_{k,N-1}} 
	\left|a_{k,j_0,\ldots,j_{N-1}}\right|^\frac{1}{k}I_{k,j_0,\ldots,j_{N-1}}(f)\right),
	\end{split}
	\end{equation*}
where
	\begin{equation*}
	\begin{split}
	&I_{k,j_0,\ldots,j_{N-1}}(f)=\prod_{\ell=0}^{N-1}\left|f^{(\ell)}\right|^{\frac{j_\ell}{k}}
	\prod_{\ell=0}^{N-2}\left(\left(f^{(\ell)}\right)^{M_\ell+1}\right)^\#
	\frac{|f^{(N-1)}|^{M_{N-1}}}{1+|f^{(N-1)}|^{2(M_{N-1}+1)}}.
	\end{split}
	\end{equation*}
Therefore to deduce the assertion it suffices to show that $I_{k,j_0,\ldots,j_{N-1}}(f)\lesssim1$. But a direct calculation shows that
	\begin{equation*}
	\begin{split}
	I_{k,j_0,\ldots,j_{N-1}}(f)&=\prod_{\ell=1}^{N-1}\frac{(M_\ell+1)|f^{(\ell)}|^{\frac{j_\ell}{k}+M_\ell+1}}{1+|f^{(\ell)}|^{2(M_\ell+1)}}
	\frac{|f|^{\frac{j_0}{k}+M_0}}{1+|f|^{2(M_0+1)}}\\
	&\le\prod_{\ell=1}^{N-1}(M_\ell+1)\frac{\left(1+|f^{(\ell)}|^{2(M_\ell+1)}\right)^{\frac{^{\frac{m_{k,\ell}}{k}+M_\ell+1}}{2(M_\ell+1)}}}{1+|f^{(\ell)}|^{2(M_\ell+1)}}
	\frac{\left(1+|f|^{2(M_0+1)}\right)^{\frac{\frac{m_{k,0}}{k}+M_0}{2(M_0+1)}}}{1+|f|^{2(M_0+1)}}\\
	&\le\prod_{\ell=1}^{N-1}(M_\ell+1),\quad k=1,\ldots,n,
	\end{split}
	\end{equation*}
by the hypotheses \eqref{hypotheses}.\hfill$\Box$

\section{Proof of Theorem~\ref{theorem:f-m-implications}}

We prove first (i). Assume on the contrary that $f\not\in\N$. Then for each $Z\in f(\D)$, with at most four possible exceptions, 
$\sup\{f^\#(z)(1-|z|^2):z\in\D,\,f(z)=Z\}=\infty$ by \cite{Lappan1974}. Let $Z\in f(\D)\setminus\{0,\infty\}$ be one of the points for which the supremum is infinity, and let $\{z_n\}_{n=1}^\infty$ denote a sequence of preimages of $Z$ such that $f^\#(z_n)(1-|z_n|^2)\to\infty$, as $n\to\infty$. Then
	\begin{equation}
	\begin{split}
	\left(f^m\right)^\#(z_n)(1-|z_n|^2)
	&=\frac{m|f(z_n)|^{m-1}(1+|f(z_n)|^2)}{(1+|f(z_n)|^{2m})}f^\#(z_n)(1-|z_n|^2)\\
	&=\frac{m|Z|^{m-1}(1+|Z|^2)}{(1+|Z|^{2m})}f^\#(z_n)(1-|z_n|^2)\to\infty,\quad n\to\infty,
	\end{split}
	\end{equation}
and therefore $f^m\not\in\N$.
A reasoning similar to that in the case (i) with \cite[Theorem~9]{AR2011} gives (iii) and (iv). Further, \cite{AW2011} together with Lappan's proof in \cite{Lappan1974} can be used to establish an analogue of the five-point theorem for strongly normal functions, which in turn gives (ii) as above. 

To prove (v)-(vii), we use \cite[Theorem~2]{Yamashita1981}. It implies 
	\begin{equation}\label{Eq:Yamashita}
	\int_{D(0,r)}f^\#(z)^2\log\frac{r}{|z|}\,dA(z)
	\lesssim\int_{D(0,r)}\left(f^m\right)^\#(z)^2\log\frac{r}{|z|}\,dA(z)+1,\quad 0<r<1.
	\end{equation}
By letting $r\to1^-$ we deduce
	\begin{equation}\label{111111}
	\int_\D f^\#(z)^2\log\frac1{|z|}\,dA(z)
	\lesssim\int_\D \left(f^m\right)^\#(z)^2\log\frac1{|z|}\,dA(z)+1.
	\end{equation}
This together with the inequalities $1-t\le-\log t\le\frac1t(1-t)$, valid for all $0<t\le1$, yield (v).

By integrating \eqref{Eq:Yamashita} over $(0,1)$ with respect to $\om(r)r\,dr$ and applying Fubini's theorem we deduce
	\begin{equation}\label{eq:plaah}
	\int_\D f^\#(z)^2\om^\star(z)\,dA(z)\lesssim\int_\D\left(f^m\right)^\#(z)^2\om^\star(z)\,dA(z)+\|\om\|_{L^1(0,1)}
	\end{equation}
from which the assertion (vi) follows.

By applying \eqref{111111} to $f\circ\vp_a$ we deduce (vii). This completes the proof of the theorem.

\end{document}